\documentclass{article}
\usepackage[pdftex]{graphicx}
\usepackage{amsmath, amssymb}
\usepackage{type1cm}
\usepackage{authblk}
\usepackage{latexsym}
\usepackage{subfigure}
\usepackage{lscape}
\usepackage{threeparttable}

\begin{document}

\title{Enumeration of associative magic squares \\ of order 7}

\author[1]{Go Kato}
\author[1]{Shin-ichi Minato}

\affil[1]{Kyoto University}

\date{\today}

\maketitle

\begin{abstract}
An associative magic square is a magic square such that the sum of any 2 cells at symmetric positions with respect to the center is constant. The total number of associative magic squares of order 7 is enormous, and thus, it is not realistic to obtain the number by simple backtracking. As a recent result, Artem Ripatti reported the number of semi-magic squares of order 6 (the magic squares of $ 6 \times 6 $ without diagonal sum conditions) in 2018. In this research, with reference to Ripatti's method of enumerating  semi-magic squares, we have calculated the total number of associative magic squares of order 7. There are exactly 1,125,154,039,419,854,784 associative magic squares of order 7, excluding symmetric patterns.
\end{abstract}
 
\section{Introduction}
\begin{figure}
\begin{minipage}
  {0.32\columnwidth}
\centering
\subfigure[3x3 magic\label{magicexam}]{\includegraphics[width=2.5cm,height=2.5cm]{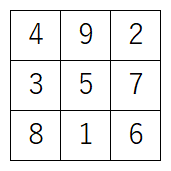}}
\end{minipage}
\begin{minipage}
  {0.32\columnwidth}
\centering
\subfigure[3x3 semi-magic]{\includegraphics[width=2.5cm,height=2.5cm]{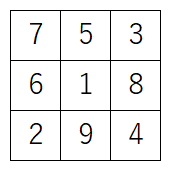}}
\end{minipage}
\begin{minipage}
  {0.32\columnwidth}
\centering
\subfigure[4x4 associative magic]{\includegraphics[width=2.5cm,height=2.5cm]{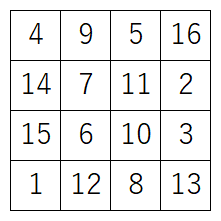}}
\end{minipage}
\label{magicexamplesfirst}
\caption{example of magic squares\label{magicexamples}}
\end{figure}

A magic square of order $n$ is an $n \times n$ square grid such that the sums of the numbers in each row, column, and diagonal are equal. The semi-magic square and associative magic square are special kinds of magic square. A semi-magic square is a magic square without the diagonal sum condition. An associative magic square is a magic square such that the sum of any 2 cells at symmetric positions with respect to the center is constant. Figure \ref{magicexamplesfirst} shows examples of these squares.

\begin{table}[tb]
\normalsize
\centering
\caption{number of magic squares}
\label{ems}
\begin{tabular}[t]{rrrr}
\hline\hline
n&semi-magic&magic&associative magic\\ \hline
3&9&1&1\\
4&68,688&880&48\\ 
5&579,043,051,200&275,305,224&48,544\\
6&94,590,660,245,399,996,601,600&(unsolved)&0\\
7&(unsolved)&(unsolved)&\bf(this research)\\ \hline
\end{tabular}
\end{table}

The known numbers of these squares are summarized in Table \ref{ems}. These squares remain of the same kind even after rotating and reflecting their entries. Thus, the table represents numbers of essentially different squares relative to these operations. The history of magic squares is long. It is said that a magic square of order 3, called lo shu, was described in China around 2200 BC. 880 magic squares of order 4 were found by Bernard Frenenicle de Bessy in 1693 \cite{4magic}, and Kathleen Ollerenshaw and Hermann Bondi proved that no other magic square of order 4 exists \cite{4magicproof}. The number of magic squares of order 5 was calculated by Richard Schroeppel by using a backtracking algorithm in 1973, and the results were published by Martin Gardner in 1976 \cite{5magic}. The number of magic squares of order 5 was calculated in 1973, but even after more than 40 years, the number of magic squares of order 6 remains unknown.

According to Walter Trump's website \cite{walter}, Mutsumi Suzuki calculated the number of associative magic squares of order 5 to be 48,544. In 1919, Charles Planck proved that there are no associative magic squares of order 6 \cite{planck}. Although it is known that there are many associative magic squares of order 7, the exact number was not known until this report.

The number of $7 \times 7$ associative magic squares was estimated to be within the range of $ (1.125151 \pm 0.000051) \times10^{18} $ with a probability of 99\% by Walter Trump, who used a method combining Monte Carlo and backtracking \cite{walter}. Since there are approximately $10^{18}$ associative squares of order 7, we cannot calculate the exact number in a realistic amount of time with simple backtracking algorithms, which take a computational time at least proportional to the number of solutions.

On the other hand, Artem Ripatti recently reported the number of semi-magic squares of order 6 in 2018 \cite{artem}. Ripatti divided the square into two parts and enumerated each part. Then, he combined the enumerations of each part. This method proved to be faster than simple backtracking.

In this paper, we extend Ripatti's method for semi-magic squares to associative magic squares and propose an algorithm to calculate the number of associative magic squares of order 7. Section \ref{pre} describes the properties of associative squares. Next, Section \ref{kizon} outlines Ripatti's method of counting $6 \times 6$ semi-magic squares. Section \ref{ours} describes our algorithm to enumerate associative magic squares of order 7, and Section \ref{resultsec} shows the results of the calculation.

\section{Preliminary\label{pre}}
\begin{figure}[b]
\centering
\includegraphics[width=4cm,height=4cm]{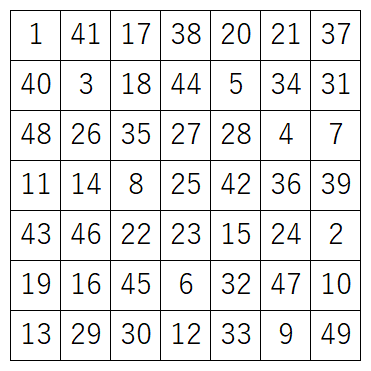}
\caption{example of an associative magic square of order 7\label{assoex}}
\end{figure}

In this paper, we call an $n \times n$ matrix such that natural numbers $ 1, 2, \cdots , n^{2} $ appear once only a square (of order $n$). Let $X_{ij}$ be the $i$-th row, $j$-th column element of the square. An associative magic square of order 7 is a magic square such that the sum of any 2 cells at symmetric positions from the center equal 50. An example of a $7 \times 7$ associative magic square is shown in Figure \ref{assoex}.
The sum of the numbers in each row of an associative magic square is $ (1 + 2 + \cdots + 49) / 7 = 175 $. Each associative magic square represents exactly 8 associative magic squares under rotation and reflection. We can calculate the total number of associative magic squares and then divide it by 8 to get the number of associative magic squares up to a reflection or rotation.

\subsection{Properties of associative squares \label{seisitu}}
A $7 \times 7$ square (not necessarily a magic square) such that the sum of any two elements at symmetrical positions from the center element is constant ($ X_{ij} + X_{(8-i) (8-j)} = 7^{2}+1=50$) has the following properties.
\begin{enumerate}
\item The central element $X_{44}$ is always 25. \label{one}
\item The element at the counterpart position can be determined as follows: \\
$ X_{ij} = 50-X_{(8-i) (8-j)} (i, j = 1, 2, \cdots, 7) $ \label{two}
\item The sum of the numbers on each diagonal must be 175. \\ 
$ \sum_{i = 1}^{7}X_{ii} = \sum_{i = 1}^{7} X_{(8-i) i} = 175$ \label{three}
\item The sums of the 4-th row and column must each be 175. \\ 
$ \sum_{i = 1}^{7} X_{i4} = \sum_{i = 1}^{7} X_{4i} = 175 $ \label{four}
\item The sum of the $k$-th row is 175 if and only if the sum of the $(8-k)$-th row is 175. \\ 
$ \sum_{i = 1}^{7} X_{ki} = 175 \Leftrightarrow \sum_{i = 1}^{7} X_{(8-k) i} = 175 (k = 1, 2 , 3) $ \label{five}
\end{enumerate}

\begin{figure}[t]
\centering
\includegraphics[width=4cm,height=4cm]{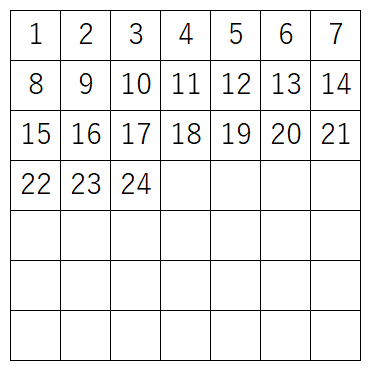}
\caption{24 key cells of associative magic squares \label{associative24}}
\end{figure}

From (\ref{one}) and (\ref{two}), all elements of the square are determined once we have decided on the values of the 24 cells shown in Figure \ref{associative24}. The sums of all the rows, columns, and diagonals should be 175 in order for a square $X$ satisfying the above symmetry to be a magic square. Moreover, for $X$ to be a magic square, it is sufficient that only 1st, 2nd, 3rd rows and columns of $X$ from (\ref{three}), (\ref{four}) and (\ref{five}) satisfy the constraints.

\section{Previous Methods \label{kizon}}
\begin{figure}[tb]
\centering
\includegraphics[width = 4 cm, height = 4 cm]{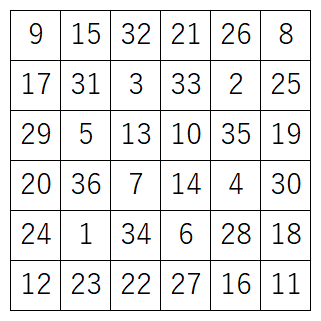}
\caption{example of a $6 \times 6$ semi-magic square \label{semiex}}
\end{figure}

In this section, we outline the method \cite{artem} that Ripatti used for counting semi-magic squares of order 6. A semi-magic square of order 6 is a $6 \times 6$ square grid such that the sums of the numbers in any row or column (not necessarily a diagonal) are equal. An example of a $6 \times 6$ semi-magic square is shown in Figure \ref{semiex}. The sum of the numbers of each row or column of a $6 \times 6$ semi-magic square is $ (1 + 2 \cdots + 36) / 6 = 111 $. The number of semi-magic squares of order 6 is 94,590,660,245,399,996,601,600, and as in the case of associative magic squares of order 7, it is not realistic to count them with backtracking.

A semi-magic square can be transformed into $6! \times 6! = 518,400$ different semi-magic squares by rearranging the rows and columns. There are $6!$ arrangements of rows and $6!$ arrangements of columns. 

Ripatti defined canonical semi-magic squares as representatives of these 518,400 semi-magic squares. We can count the canonical semi-magic squares and then multiply the number by 518,400 to get the total number of semi-magic squares.

\begin{figure}[tb]
\centering
\includegraphics[width = 4 cm, height = 4 cm]{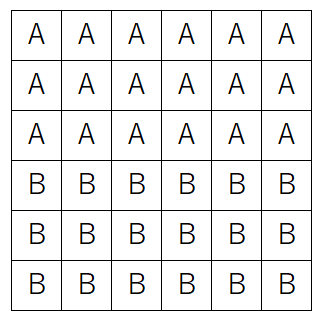}
\caption{division of a $6 \times 6$ square \label{semidiv}}
\end{figure}
\begin{figure}
\centering
\includegraphics[width = 4cm, height = 4cm]{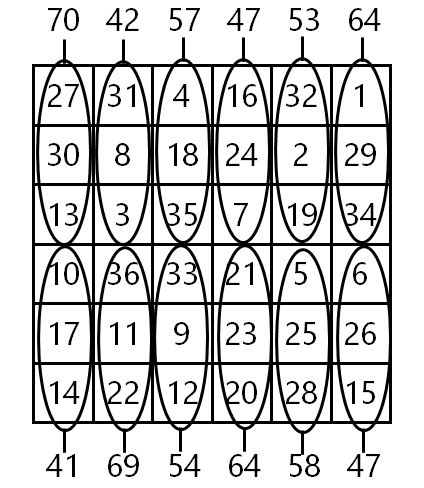}
\caption {profile of a semi-magic square \label{semiprofile}}
\end{figure}

Next, Ripatti divided the square into the upper half and lower half, as shown in Figure \ref{semidiv}. Let the $profile$ of the upper half of the square be:
\begin{eqnarray*}
(\sum_{i=1}^{3}X_{i1}, \:\sum_{i=1}^{3}X_{i2}, \:\sum_{i=1}^{3}X_{i3},\:\cdots, \:\sum_{i=1}^{3}X_{i6})
\end{eqnarray*}
Let the profile of the lower half of the square be:
\begin {eqnarray*}
(111-\sum_{i=4}^{6}X_{i1},\:111-\sum_{i=4}^{6}X_{i2},\:111-\sum_{i=4}^{6}X_{i3},\:\cdots, \:111-\sum_{i=4}^{6}X_{i6})
\end {eqnarray*}
For example, the profile of the upper half of Figure \ref{semiprofile} is $ (70, 42, 57, 47, 53, 64) $, while that of the lower half is $(111-41,111-69,111-54,111-64,111-58,111-47) = (70, 42, 57, 47, 53, 64)$. In accordance with these definitions, the combination of the upper half and lower half is a square such that the sum of the elements of each column is 111 when the two halves have the same profile. Therefore, the combination of the upper half and lower half squares is a semi-magic square if and only if the combination satisfies all of the following conditions.
\begin{itemize}
\item Each of the upper half and lower half squares satisfy the semi-magic row constraints (the sum of the elements of each row is equal to 111).
\item The profiles of the upper and lower squares are the same.
\item Each element appears only once in the combination of the upper half and lower half.
\end{itemize}
It is possible to count the number of semi-magic squares according to the following procedure.

\begin{enumerate}
\item Create all sets $(U, L)$ ($U$ is a set of 18 elements in the upper half of the square, and $L$ is a set of the 18 elements of the lower half.). Perform the following procedure for each $(U, L)$. \\
($ U \cup L = \{1, 2, \cdots, 36 \}$)
\item Prepare an array of counters $N_{U}[p]$ for a respective profile $p$, and initialize its elements to 0. Generate all the upper halves of the square such that the sum of the elements of each row is 111 using $U$ numbers; then compute the profile value $p$ and increment $N_{U}[p]$ for each of the upper half. Thus, we can calculate that there are $N_{U}[p]$ patterns for the upper half satisfying the constraints for each profile. \label{semicounttwo}
\item As well as (\ref{semicounttwo}), prepare an array of counters $N_{L}[p]$ for a respective profile $p$, and initialize its elements to 0. Generate all the lower halves of the square such that the sum of each row is 111 using $L$ numbers; then compute the profile value $p$ and increment $N_{L}[p]$ for each of the lower half. \label{semicountthree}
\item The number of semi-magic squares of order 6 for one pair $(U, L)$ is \\$ \sum_{p} (N_{U}[p] \times N_{L}[p]) $. \label{semicountfour}
\end{enumerate}

In this way, we can divide the problem of enumerating semi-magic squares into two problems: counting the upper half of the squares and counting the lower half. If we use simple backtracking without splitting the square, we need to count $N_{U}[p] \times N_{L}[p]$ semi-magic squares one by one. However, when counting separately, we count only $N_{U}[p]$ upper halves and $N_{L}[p]$ lower halves in order to count $N_{U}[p] \times N_{L}[p]$ semi-magic squares.

Straightforwardly, the number of pairs $(U, L)$ is ${}_{36} C _{18}$, but by considering constraints such as $\sum_{x \in U}x = \sum_{x \in L}x = 111 \times 3$, Ripatti suggested that it is sufficient to consider only 9,366,138 $(U, L)$ pairs.

According to Ripatti, we can efficiently count the upper halves of the squares and lower halves by using the family of sets of numbers that fit on individual lines of semi-magic squares, \mbox{\boldmath $S$}$\:\:=\{s\:|\:s\subset \{1,2,\cdots,36\}, \:\:|s|=6, \:\:\sum_{a \in s}a=111\}$.

It is possible to calculate efficiently if we count only the combinations of the upper half square and the lower half square which become canonical semi-magic squares. Ripatti also used other speed-up techniques and handled corner cases, but we will omit them here.

Speeding up the calculations by parallelization is effective since this counting method is completely independent for each $(U, L)$ pair. Ripatti'{}s method took over 5 months on 10 threads to report that the number of semi-magic squares of order 6 was 94,590,660,245,399,996,601,600.

\section{Our method \label{ours}}
\subsection {Overview}
We extend the existing method described in Section \ref{kizon} so that it can be used to count associative magic squares of order 7. First, we describe row and column rearrangements that transform associative magic squares into different ones. Next, we propose one division of the $7 \times 7$ square to count parts separately. We define a profile well suited for the division and show that we can divide up the problem of counting associative magic squares into two. Finally, we present the detailed calculation procedure for counting associative magic squares of order 7.

\subsection {Associative magic square transformations \label {trans}}
\begin{figure}[tb]
\centering
\includegraphics[width = 8cm]{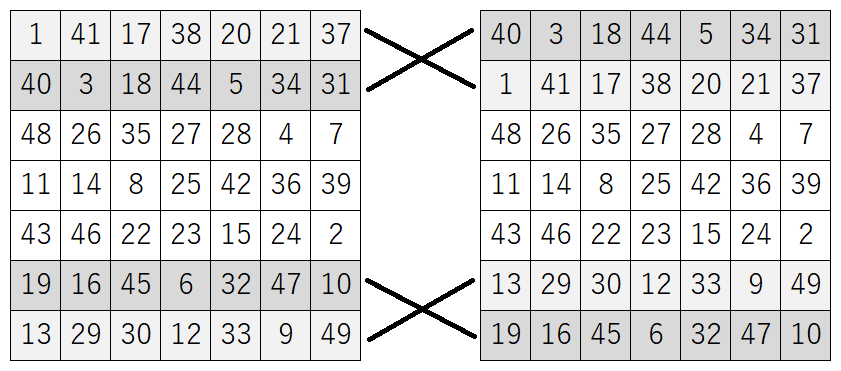}
\caption {Swapping rows (\ref {swap1267}) of an associative magic square \label{7x7swapex}}
\end{figure}

Associative magic squares of order 7 can be transformed into other associative magic squares by symmetrical swapping of rows and columns with respect to the center. Such swappings are shown below.
\begin{enumerate}
  \item Swap the 1st and 7th rows. \label{swap17}
  \item Swap the 2nd and 6th rows. \label{swap26}
  \item Swap the 3rd and 5th rows. \label {swap 35}
  \item Swap the 1st and 2nd rows, and swap 6th and 7th rows. \label{swap1267}
  \item Swap the 2nd and 3rd rows, and swap 5th and 6th rows. \label{swap2356}
  \item Swap the 1st and 7th columns.
  \item Swap the 2nd and 6th columns.
  \item Swap the 3rd and 5th columns.
  \item Swap the 1st and 2nd columns, and swap 6th and 7th columns.
  \item Swap the 2nd and 3rd columns, and swap 5th and 6th columns.
\end{enumerate}
For example, as shown in Figure \ref{7x7swapex}, the associative magic square on the right is obtained by swapping rows (\ref {swap1267}) of the left square.

Since the above rearrangements are symmetrical replacements of rows and columns, two elements in the symmetrical position before swapping remain symmetrical after swapping. As stated in (\ref{three}) in Section \ref{seisitu}, the sums of the numbers on the diagonals equal 175 when the sum of any two symmetrical elements is constant. The sum of each row and each column after transformation is the same as before, since the set of seven numbers in each row and each column remains the same before and after the replacements. Therefore, associative magic squares remain associative after symmetrical rearrangements of rows and columns such as the above.

\begin{figure}[p]
\centering
\includegraphics[width = 11cm]{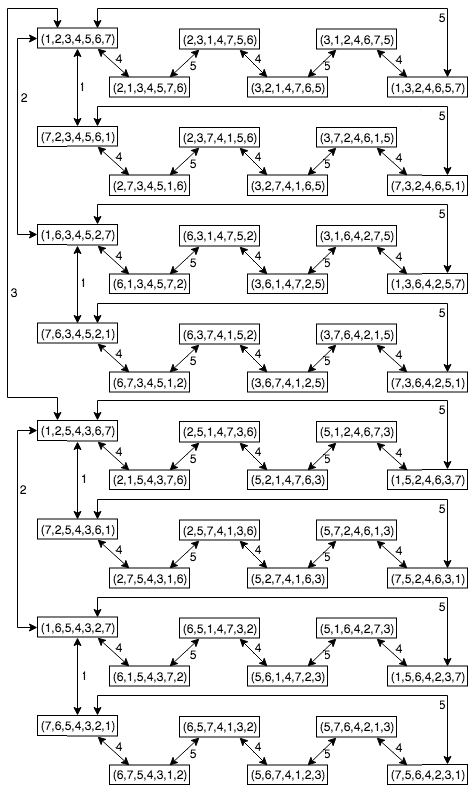}
\caption {48 permutations of rows in associative magic squares of order 7 \label{rowswapall}}
\end{figure}

Figure \ref{rowswapall} shows how symmetric squares are transformed by the above row swapping. $(r_{1}, r_{2}, r_{3}, \cdots, r_{7})$ represents an associative magic square such that the $i$-th row is replaced by the $r_{i}$-th row of the original associative magic square. The original associative magic square is represented by $(1,2,3,4,5,6,7)$. \\$(r_{1}, r_{2}, r_{3}, \cdots, r_{7}) = (r_{1}, r_{2}, r_{3}, 4, 8-r_{3}, 8-r_{2}, 8-r_{1})$, because the two rows in the symmetrical position remain symmetrical with row swapping (\ref{swap17}) 〜 (\ref{swap2356}), and the fourth row doesn't move. Thus, three pairs $ (1,7), (2,6), (3,5) $ are assigned to one of $ r_ {1}, r_ {2}, r_ {3} $, respectively. There are $ 2 ^ {3} $ ways to choose which of the two numbers of each pair is assigned to $ r_{1}, r_{2}, r_{3} $, which is expressed as the vertical transformations shown in Figure \ref{rowswapall}. There are $3!$ ways to choose pairs to be assigned to one of $(r_{1}, r_{7}), (r_{2}, r_{6}), (r_{3}, r_{5})$, which is expressed as the horizontal transformations shown in Figure \ref{rowswapall}. Therefore, an associative magic square can be transformed $2^{3} \times 3! = 48$ ways by rearranging the rows. We can also transform an associative magic square into 48 associative magic squares by rearranging the columns. Thus, an associative magic square of order 7 can be transformed into $48 \times 48 = 2304$ associative magic squares in total.

In Section \ref{hyoujunkei}, we defined a canonical associative magic square of order 7 that represents the 2304 associative magic squares made by such transformations. We can count all the associative magic squares by calculating the number of only canonical associative magic squares and multiplying that number by 2304.

\subsection{Square division}
Since there are 49 cells in the $7 \times 7$ square, there are $ 2 ^ {49} $ ways to divide the square into two. Some examples are shown below.

\begin{figure}
\begin{minipage}
  {0.24\columnwidth}
\centering
\subfigure[3 upper rows and 4 lower rows \label{div1}]
{\includegraphics[width = 2.5cm, height = 2.5cm]{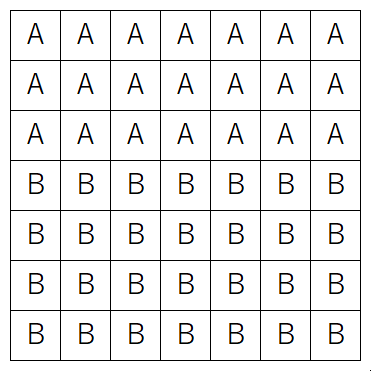}}
\end{minipage}
\begin{minipage}
  {0.24\columnwidth}
\centering
\subfigure[checkered pattern\label{div4}]
{\includegraphics[width = 2.5cm, height = 2.5cm]{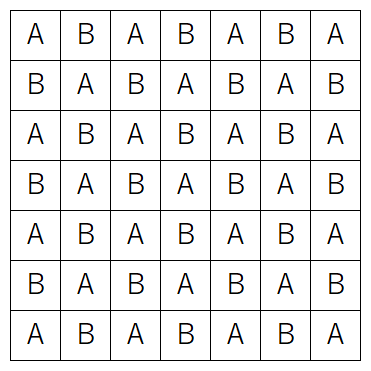}}
\end{minipage}
\begin{minipage}
  {0.24\columnwidth}
\centering
\subfigure[1 center row and 6 outer rows \label{div2}]
{\includegraphics[width=2.5cm,height=2.5cm]{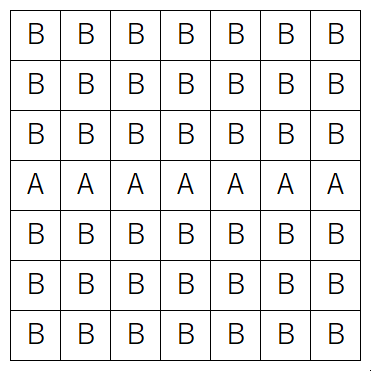}}
\end{minipage}
\begin{minipage}
  {0.24\columnwidth}
\centering
\subfigure[3 center rows and 4 outer rows \label {div3}]
{\includegraphics[width=2.5cm,height=2.5cm]{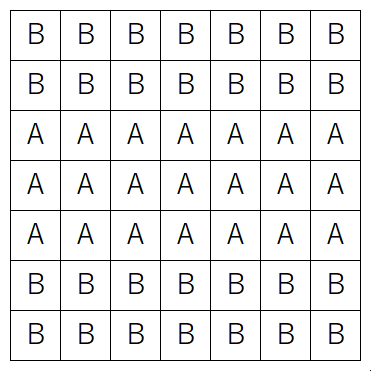}}
\end{minipage}
\label{magicdivisions}
\caption{examples of dividing the square}
\end{figure}

If we divide the square into an upper half and lower half like the semi-magic square division shown in Figure \ref{div1} to count associative magic squares, it is difficult to consider the upper half and lower half squares independently because the arrangement of elements in the upper half complexly affects the arrangement of elements in the lower half, wherein $ X_{ij} + X_{(8-i) (8-j)} = 50 $. Therefore, it is preferable to divide up a square so that two cells at symmetrical positions are included in the same group.

If we divide the square into two alternating groups, A and B, as shown in Figure \ref{div4}, two cells at symmetrical positions are put into the same group. However, each row and each column of the square are divided into two groups. It is difficult to consider groups A and B independently because the arrangement of elements in group A complexly affects the arrangement of elements in group B such that the sum of the elements of each row and column is equal to $175$. Thus, it is better not to divide the rows and columns in the two groups as much as possible.

As shown in Figure \ref{div2}, when dividing the square into one center row and six outer rows, symmetrical cells are included in the same group, and the rows are the not divided. However, the number of cells in the six outer rows is too large. Therefore, there is little difference in problem size between counting associative squares of order 7 directly and counting 6 outer rows that satisfy some constraints. As a result, we cannot efficiently calculate with such divisions. We should divide the square almost into two halves in order to reduce the sizes of the two counting subproblems obtained by the division.

From the above, the division shown in Figure \ref{div3} is superior to many of the other divisions. Therefore, we will consider counting associative magic squares of order 7 by dividing up the square into 3 center rows and 4 outer rows.

\begin{figure}[b]
\begin{minipage}
  {0.49\columnwidth}
\centering
\subfigure[profile of the center part \label{asoprofile1}]
{\includegraphics[width=3.5cm,height=3.5cm]{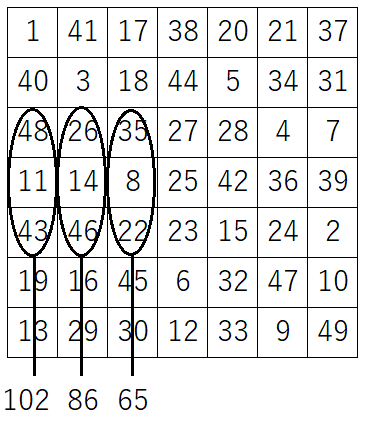}}
\end{minipage}
\begin{minipage}
  {0.49\columnwidth}
\centering
\subfigure[profile of the outer part\label{asoprofile2}]
{\includegraphics[width=3.5cm,height=3.5cm]{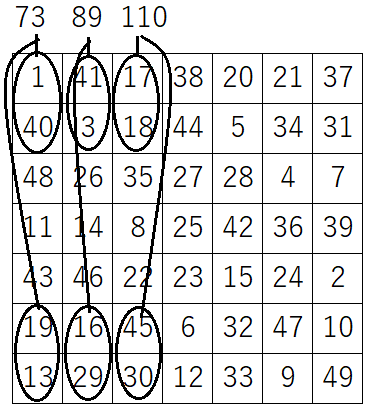}}
\end{minipage}
\caption{profiles of $7 \times 7$ associative magic square}
\end{figure}

\subsection{Profiles of divided groups \label{gattai}}
Considering the arrangements of the numbers in the center part and outer part such that the sums of any two symmetrical cells are the same, from the considerations given in Section \ref{seisitu}, a combination of arrangements of the center part and the outer part becomes an associative magic square if and only if sums of the elements of the individual 1st, 2nd, and 3rd rows and columns are 175, and each element appears only once in the combination. In the existing method of the semi-magic squares of order 6, it was necessary to define the profile as a vector of 6 partial sums of the columns, but we can define the profiles of associative magic squares as a vector of 3 partial sums for the left three columns because of their symmetric constraints. We define the profile of the center parts as 
\begin{eqnarray*}
(\sum_{i = 3}^{5} X_{i1}, \: \: \sum_{i = 3}^{5}X_{i2}, \: \: \sum_{i = 3}^{5}X_{i3})
\end{eqnarray*}
and define profile of the outer parts as
\begin{eqnarray*}
(175-X_{11} + X_{21} + X_{61} + X_{71}, \\ 175-X_{12} + X_{22} + X_{62} + X_{72}, \\ 175 -X_{13} + X_{23} + X_{63} + X_{73})
\end{eqnarray*}
For example, the profile of the center part of the square shown in Figure \ref{asoprofile1} is $ (102,86,65) $, and the profile of the outer part of the same square (Figure \ref{asoprofile2}) is $(175-73,175-89,175-110)$. The sums of the 1st, 2nd and 3rd columns of a combination of center and outer rows are each 175 if and only if the combination consists of center and outer parts that have the same profile.

Therefore, a combination of the center part and outer part of a square is an associative magic square if and only if it satisfies all of the following conditions.
\begin{itemize}
\item For each of the center part and outer part of the square, the sum of any 2 symmetrical elements is constant.
\item The sum of the elements of the 3rd row of the square is equal to 175.
\item The sums of the elements of the 1st and 2nd rows of the square are each equal to 175.
\item The center part and the outer part of the square have the same profile.
\item Each element appears exactly once in the combination.
\end{itemize}

\subsection {Procedure for counting associative magic squares of order 7 \label{count}}
From the discussion in Section \ref{gattai}, we can calculate the total number of associative magic squares by using the following procedure.
\begin{enumerate}
\item Create all sets $(S_A, S_B)$ ($S_A$ is a set of 21 elements in the center parts of the square, and $S_B$ is a set of 28 elements in the outer parts.). Perform the following procedure for each $(S_A, S_B)$. \\ ($ S_A \cup S_B = \{1, 2, \cdots, 49 \} $) \label{countone}
\item Prepare an array of counters $N_A[p]$ for the respective profile $p$, and initialize its entries to 0. Using $S_A$ numbers, generate all the center parts of the square such that the sum of the elements of the 3rd row of the square is 175 and sum of any two symmetrical elements is 50, satisfying canonical associative magic square constraints. For each of the generated center parts of the square, compute the profile value $p$ and increment $N_{A}[p]$. Thus, we can calculate that there are $N_{A}[p]$ patterns for the center part satisfying the constraints for each profile. \label{counttwo}
\item As well as (\ref{counttwo}), prepare an array of counters $N_{B}[p]$ for respective profile $p$, and initialize its entries to 0. Using $S_B$ numbers, generate all the outer parts of the square such that the sums of the elements of the 1st and 2nd rows are each 175 and the sum of any two symmetrical elements is 50, satisfying the canonical associative magic square constraints. For each of the outer parts, compute the profile value $p$ and increment $N_{B}[p]$ . \label{countthree}
\item The number of canonical associative magic squares of order 7 for one pair $ (S_A, S_B) $ is $ \sum_{p}(N_{A}[p] \times N_{B}[p])$. \label{countfour}
\end{enumerate}

The details of (\ref{countone}) are in Section \ref{countdivce}, while the definitions of canonical associative magic squares are in Section \ref{hyoujunkei}. Moreover, the details of (\ref{counttwo}) are in Section \ref{counttyuuou}, and the details of (\ref{countthree}) are described in Section \ref{counthasi}.

In this way, we can divide the problem of enumerating associative magic squares of order 7 into two problems: counting the center parts of the square and counting the outer parts. If we enumerate them with simple backtracking without splitting the square, we need to explore $N_{A}[p] \times N_{B}[p]$ patterns of associative magic squares one by one. However, when counting halves separately, we explore only $N_A[p] + N_B[p]$ patterns to count $N_{A}[p] \times N_{B}[p]$ associative magic squares.

\subsubsection {Number of $(S_A, S_B)$ pairs \label{countdivce}}
We will discuss the number of $(S_A, S_B)$ pairs in the procedure of (\ref{countone}) in Section \ref{count}. Considering simply, the number of pairs $ (S_A, S_B) $ is $ {}_{49}C_{21} $ (about $ 3.9 \times 10^{13} $). But, there are $ {}_{48}C_{20} $ ways, if we consider that $X_{44}$ must be $25$. Since the sum of any two symmetrical elements about the center is constant, if $ S_A $ includes $ x $, then $ S_A $ must include $ 50-x $. Therefore, it is sufficient to divide 24 numbers $ 1, 2, \cdots, 24 $ into $ (S_A, S_B) $. Thus, we only have to calculate ${}_{24}C_{10} = 1, 961, 256 $ $ (S_A, S_B) $ pairs.

\subsubsection {Canonical associative magic squares of order 7 \label{hyoujunkei}}
As described in Section \ref{trans}, we define a canonical associative magic square of order 7 that represents the 2304 associative magic squares made by the transformations, and count only the canonical associative magic squares. We can calculate the number of associative magic squares by multiplying the number of canonical squares by 2304 and thus reduce the number of squares that need to be counted by a factor of 2304.

In the procedure shown in Section \ref{count}, it takes much more computation time to count the outer parts than the center parts because of the difference in the number of included cells. Therefore, we face a bottleneck in counting the outer parts that satisfy some conditions. We can prune the search to count the outer parts efficiently by imposing the constraints of canonical associative magic squares as much as possible on the outer parts. From the above, we define canonical associative magic squares as those that satisfy all of the following conditions.

\begin{enumerate}
\item $X_{11} < X_{17}, X_{12} < X_{16}, X_{13} < X_{15}$\label{transone}
\item $X_{11} < X_{12} < X_{13}$\label{transtwo}
\item $row_{1} < row_{7}, row_{2} < row_{6}, row_{3} < row_{5}$\label{transthree}
\item $row_{3} < row_{1} < row_{2}$\label{transfour}
\end{enumerate}
However, $row_{i}=\min\{X_{i1},X_{i2},\cdots,X_{i7}\}$.

(\ref{transone}) and (\ref{transtwo}) are restrictions on column swapping, and these restrictions constrain the first row. (\ref {transthree}) and (\ref{transfour}) are restrictions on row swapping, but it is not possible to put all these restrictions on the four outer rows. We make $row_{3}$ small instead, and the following formula holds:
\begin{eqnarray*}
row_{3} = \min\{row_{1}, row_{2}, row_{3}, row_{5}, row_{6}, row_{7}\}
\end {eqnarray*}
The 3rd row of the canonical squares always contains the smallest number, except for numbers in the 4th row. Thus, the center parts always contain the element one. The number of $ (S_A, S_B) $ pairs that need to be calculated, as described in Section \ref{countdivce}, is found to be $ {}_{23}C_{9} = 817, 190 $ using this property of canonical squares. By defining the canonical associative magic squares in this way, the counting of the outer parts becomes approximately 384 times faster, and the number of $ (S_A, S_B) $ pairs that needs to be calculated is approximately 0.4 times the total.

\subsubsection {Counting the center parts \label{counttyuuou}}
Here, we describe the (\ref{counttwo}) of the procedure in Section \ref{count}. From Figure \ref{associative24} in Section \ref{seisitu}, if the numbers in the 3rd row of the square and $ X_{41}, X_{42}, X_{43} $ are determined, the remaining elements in the center part are also determined.

We define \mbox{\boldmath $R$}, a family of sets of 7 numbers that can be put in rows other than the 4th row of the associative magic of order 7, as
\begin {eqnarray*}
\mbox{\boldmath $R$} = \{r \: | \: r \subset \{1,2, \cdots, 24,26,27, \cdots, 49 \}, \: | r | \: = 7, \\ \sum_{ a \in r} a = 175, \forall x, y \in r, x + y \neq 50 \}
\end {eqnarray*}
We have generated \mbox{\boldmath $R$} in preparation for counting the center parts. The family size of \mbox{\boldmath $R$} is 452,188. When we use $S_A$ numbers, \mbox{\boldmath $R_3$}, the family of possible sets of 7 numbers in the 3rd row of an associative square of order 7, can be easily calculated as 
\begin{eqnarray*}
\mbox{\boldmath $R_3$}= \{r \: | \: r \in \mbox{\boldmath $R$}, r \subset S_A \}
\end{eqnarray*}

For each $r \in \mbox{\boldmath $R_3$}$, there are $7!$ ways for all permutations of $r$ to choose the elements of the 3rd row, and the elements in the 5th row are determined by the elements in the symmetrical 3rd row. But, we exclude the $ row_{3} < row_{5} $ arrangement from the constraints of the canonical form. There are $6 \times 4 \times 2$ ways to choose the elements of $X_{41},X_{42},X_{43}$ from the unused numbers in $S_A$. We generate all the arrangements of elements in the center parts, calculate of profile $ p $ for each arrangement, and increment $ N_{A}[p] $.

\subsubsection{Counting the outer parts \label{counthasi}}
Here, we describe (\ref{countthree}) of the procedure in Section \ref{count}. From Figure \ref{associative24} in Section \ref{seisitu}, if the elements in the 1st and 2nd rows are determined, the remaining elements in the outer parts are also determined.

Besides counting the center parts, \mbox{\boldmath $R_1$}, we can easily calculate the family of sets of 7 elements in the first row of the associative magic square, 
\begin{eqnarray*}
\mbox{\boldmath $R_1$}= \{r \: | \: r \in \mbox{\boldmath $R$}, r \subset S_B \}
\end{eqnarray*}
Next, we define \mbox{\boldmath $R_{26}$}, a family of 14-number subsets of $S_B$ excluding the numbers in the 1st and 7th rows, as
\begin{eqnarray*}
\mbox{\boldmath $R_{26}$} = \{r \: | \: r \subset \{1, 2, \cdots, 24, 26, 27, \cdots, 49 \}, \: | r | = 14, \\ \forall x \in r, \: \exists y \in r, \: x + y = 50 \}
\end{eqnarray*}
Furthermore, we can pre-generate \mbox{\boldmath $R_2$} for all \mbox{\boldmath $R_{26}$}, a family of 7-number subsets of \mbox{\boldmath $R_{26}$} that can be put on the 2nd row of the associative magic square of order 7, defined as 
\begin{eqnarray*}
\mbox{\boldmath $R_2$} = \{r \: | \: r \subset \mbox{\boldmath $R_{26}$}, \: \: \: | r | = 7, \: \: \: \sum_{a \in r} a = 175, \: \: \: \forall x, y \in r, \: \: \: x + y \neq50 \}
\end{eqnarray*}
By pre-generating \mbox{\boldmath $R_2$} for all \mbox{\boldmath $R_{26}$}, it is possible to count the outer parts by the following procedure. For each $r \in \mbox{\boldmath $R_1$}$, there are $7!$ ways for all permutations of $r$ to choose the elements of the 1st row, and the elements in the 7th row are determined by the elements in the symmetrical 1st row. But, we exclude arrangements violating the constraints of the canonical form. Let \mbox{\boldmath $R_{26}$} be the set of unused numbers of $ S_B $ (excluded $ r $ and the numbers in the 7th row). For each $r' \in \mbox{\boldmath $R_2$}$ calculated from \mbox{\boldmath $R_{26}$}, there are $7!$ ways for all permutations of $r'$ to choose the elements of the 2nd row, and the elements in the 6th row are determined by the elements in the symmetrical 2nd row. We generate all the arrangements of elements in the outer parts, calculate the profile $ p $ for each arrangement, and increment $ N_{B}[p] $.

\section{Experimental results\label{resultsec}}
\subsection {Results}
We wrote a \verb!C++! program implementing the procedure described in section \ref{ours}. We assigned $ 1, 2, \cdots, 817190 $ ID numbers to each pair of $(S_A, S_B)$. Our program counts associative magic squares of order 7 for each $(S_A, S_B)$ designated with ID numbers.

As the experimental environment, two computers were used: a PC with an Intel Core i7-4960X 3.6GHz CPU and 64GB RAM running 64bit Windows 7, and a Mac with an Intel Core i5 1.8GHz CPU and 8GB RAM running 64bit macOS High Sierra.

The calculation can be speed up by parallelization because the counting method is completely independent for each $(S_A, S_B)$ pair. We divided 817,190 $ (S_A, S_B) $ pairs into 16 groups and calculated all associative magic squares of order 7 of one group in one thread. The calculation was executed on 12 threads on the Windows computer and 4 threads on the Mac computer. The calculation took about 2 weeks, and the results are shown in Table \ref{result}.

\begin{table}[t]
\normalsize
\centering
\caption{calculation results}
\label{result}
\begin{tabular}[t]{rrrr}
\hline\hline
$(S_A, S_B)$ ID&number of associative magic&calc. time & computer\\ \hline
1-50000&100798108317305280&14.5 days&Mac\\
50001-100000&91535720218951104&14.4 days&Windows\\
100001-150000&88372685889123552&14.2 days&Windows\\
150001-200000&83733351186221856&14.1 days&Windows\\
200001-250000&81588443264793504&14.0 days&Mac\\
250001-300000&79361704382078592&13.8 days&Windows\\
300001-350000&68614934779440864&13.9 days&Windows\\
350001-400000&60333826371280992&13.9 days&Windows\\
400001-450000&58972609900819872&13.6 days&Mac\\
450001-500000&56989665917916192&13.6 days&Windows\\
500001-550000&58076327642080032&13.6 days&Windows\\
550001-600000&58605580160376480&13.7 days&Windows\\
600001-650000&56406391669618560&12.9 days&Mac\\
650001-700000&56103389221682304&13.5 days&Windows\\
700001-750000&54683346217110336&13.0 days&Windows\\
750001-817190&70977954281055264&14.8 days&Windows\\ \hline
\end{tabular}
\end{table}

We found that the number of associative magic squares of order 7 is \\1,125,154,039,419,854,784 up to a reflection or rotation.

We did not think so deeply about how to divide the 817,190 problems into 16 groups, but we did achieve an effective distribution, because the maximum calculation time was 14.8 days and the minimum calculation time was 12.9 days.

Additionally, we submitted our results to the On-Line Encyclopedia of Integer Sequences (OEIS), and it was accepted on December 10, 2018. It is currently on the website \cite{oeis}.

\subsection {Verification}
Walter Trump \cite{walter} estimated the number of associative magic squares of order 7 to be within the range $ (1.125151 \pm 0.000051) \times 10^{18} $ with a probability of 99 \%. Our result, 1,125,154,039,419,854,784, is within the range of this estimate.

Walter Trump confirmed the number of $7 \times 7$ associative magic squares of one $(S_A, S_B)$ with backtracking. Trump also confirmed our results with his own program based on our method.

As described in Section \ref{trans}, a certain associative magic square can be transformed into 2304 associative magic by swapping rows and columns. When we add a 90-degree rotation to it, a certain associative magic can be transformed to $ 2304 \times 2 = 4608 $ associative magic. Therefore, the number of associative magic squares of order 7 up to a reflection or rotation is a multiple of $ 576 (= 4608/8) $. Our result also has this property.

Our results have been confirmed by a probabilistic estimation, Trump, and the properties of associative magic squares.

\section{Concluding Remarks}
We proposed a method to count the total number of associative magic squares of order 7 by extending Ripatti's method of counting $6 \times 6$ semi-magic squares to $7 \times 7$ associative squares. The proposed method divides the $7 \times 7$ square into two and divides the problem into two smaller ones. It is important to divide the $7 \times 7$ matrix into a center part and outer part and consider them independently. The proposed method counts only canonical associative magic squares which are representative of 2304 other associative magic squares and these canonical associative magic squares depend on how the square is divided up. Our calculation shows that the total number of associative magic squares of order 7 is 1,125,154,039,419,854,784, up to a reflection or rotation.

\begin{landscape}
\begin{table}[t]
\large
\centering
\caption{known numbers of squares including the results of this research}
\label{final}
\begin{threeparttable}
\begin{tabular}[t]{rrrr}
\hline\hline
n&semi magic&magic&associative magic\\ \hline
3&9&1&1\\
4&68,688&880&48\\ 
5&579,043,051,200&275,305,224&48,544\\
6&94,590,660,245,399,996,601,600&$(1.775399(42)\times 10^{19})$&0\\
7&$(4.2848 (17)\times 10^{38})$&$(3.79809 (50) \times 10^{34})$&1,125,154,039,419,854,784\\
8&$(1.0806 (12) \times 10^{59})$&$(5.2225 (18) \times 10^{54})$&$(2.5228 (14) \times 10^{27})$\\
9&$(2.9008 (22) \times 10^{84})$&$(7.8448 (38) \times 10^{79})$&$(7.28 (15) \times 10^{40})$\\
10&$(1.4626 (16) \times 10^{115})$&$(2.4149 (12) \times 10^{110})$&$0$\\ \hline
\end{tabular}
\begin{tablenotes}\footnotesize
\item[*] $ 1.775399(42) \times 10^{19} $ and other estimated values mean that the exact number of the squares is within range $(1.775399\pm0.000042)\times 10^{19}$ with a probability of 99\%.
\end{tablenotes}
\end{threeparttable}
\end{table}
\end{landscape}

Table \ref{final} summarizes results on known squares up to $ n = 10, $ including our own. Charles Planck proved that there is no associative magic square of order $N$ ($N$ is even number that is not a multiple of 4) \cite{planck}. The table includes number of squares estimated by Walter Trump for squares whose exact number is unknown \cite{walter}.

In the future, we could try calculating the number of higher-order associative magic squares and other magic squares whose total number is unknown. However, the number of associative magic squares of higher order would be too large to count with our method. The exact number of $6 \times 6$ magic squares is very interesting, but it seems to be difficult to count by making a simple extension of our method. We may also consider other kinds of number assignment or constraint satisfaction problems.

\section*{Acknowledgements}
We would like to thank Walter Trump for personal communications confirming our results. We also thank Fran\c{c}ois Le Gall and Suguru Tamaki for their helpful comments.

\end{document}